\documentclass[11pt]{pjm}
\usepackage[english]{babel}
\usepackage{amsmath}
\usepackage[T1]{fontenc}

\newtheorem {theorem} {Theorem}
\newtheorem {proposition} {Proposition}
\newtheorem {example} {Example}
\newtheorem {definition} {Definition}
\newtheorem {lemma} {Lemma}
\newtheorem {corollary}{Corollary}
\newtheorem {remark} {Remark}

\def\un{{\mathrm{1\kern -0.30em I}}}
\def\R{\mathrm{I\hspace{-0.20em}R}}
\def\N{\mathrm{I\hspace{-0.20em}N}}
\def\K{\mathrm{I\hspace{-0.20em}K}}
\def\Z{{\mathsf{Z}\hspace{-0.30em}\mathsf{Z}}}
\def\C{\mathrm{C}\hspace{-0.50em}\mathsf{l}\,\,}

\newcommand{\bpp}{\begin{proposition}}
\newcommand{\epp}{\end{proposition}}
\newcommand{\bdf}{\begin{definition}}
\newcommand{\edf}{\end{definition}}
\newcommand{\brm}{\begin{remark}}
\newcommand{\erm}{\end{remark}}
\newcommand{\bxm}{\begin{example}}
\newcommand{\exm}{\end{example}}
\newcommand{\bcl}{\begin{corollary}}
\newcommand{\ecl}{\end{corollary}}
\newcommand{\blm}{\begin{lemma}}
\newcommand{\elm}{\end{lemma}}

\newcommand{\abs}[1]{\left|#1\right|}

\newcommand{\ind}[1]{\mathrm{1\kern -0.28em I}_{#1}}

\newcommand{\Hun}{H^1\left(\Omega\right)}

\newcommand{\Lp}[1]{L^{#1}\left(\Omega\right)}

\newcommand{\Om}{\Omega}

\newcommand{\norm}[1]{\left\|#1\right\|}

\newcommand{\veps}{\varepsilon}

\newcommand{\dsp}{\displaystyle}

\newcommand{\dvg}{\operatorname{div}}

\newcommand{\grad}{\nabla}

\newcommand{\bproof}{{\bf Proof.}\hspace{0.1cm}}
\newcommand{\eproof}{\vspace{0.3cm}}

\newcommand{\ben}{\begin{enumerate}}
\newcommand{\een}{\end{enumerate}}
\newcommand{\bq}{\begin{equation}}
\newcommand{\eq}{\end{equation}}
\newcommand{\bqs}{\begin{equation*}}
\newcommand{\eqs}{\end{equation*}}
\newcommand{\bqns}{\begin{eqnarray*}}
\newcommand{\eqns}{\end{eqnarray*}}
\newcommand{\bqn}{\begin{eqnarray}}
\newcommand{\eqn}{\end{eqnarray}}
\newcommand{\bit}{\begin{itemize}}
\newcommand{\eit}{\end{itemize}}
\newcommand{\bde}{\begin{description}}
\newcommand{\ede}{\end{description}}

\newcommand{\db}[2]{\left[\kern -0.08em\left[#1,#2\right]\kern -0.1em\right]}

\begin{document}

\title[generalized Sobolev algebras and applications]{On generalized Sobolev algebras and their applications}


\author[S. Bernard]{S\'everine Bernard}
\address{Universit\'e des Antilles et de la Guyane, Laboratoire
AOC, Campus de Fouillole, 97159 Pointe-\`a-Pitre cedex,
Guadeloupe.} \email{severine.bernard@univ-ag.fr}
\urladdr{www.univ-ag.fr/aoc}
\author[S.P. Nuiro]{Silv\`ere P. Nuiro}
\address{Universit\'e des Antilles et de la Guyane, Laboratoire
GRIMAAG, Campus de Fouillole, 97159 Pointe-\`a-Pitre cedex,
Guadeloupe.} \email{paul.nuiro@univ-ag.fr}
\urladdr{www.univ-ag.fr/grimaag}


\thanks{The authors are pleased to acknowledge A.
Delcroix for helpful advice and suggestions.}


\date{Received date / Revised version date}

\begin{abstract}

In the last two decades, many algebras of generalized functions
have been constructed, particularly the so-called generalized
Sobolev algebras.  Our goal is to study the latter and some of
their main properties. In this framework, we pose and solve a
nonlinear degenerated Dirichlet problem with irregular data such
as Dirac generalized functions.
\end{abstract}

\maketitle

{\bf{Key words }}: nonlinear degenerate Dirichlet problem,
generalized solution, Sobolev algebra, non positive solution.\\

{\bf{2000 MSC }}: 35J70, 46F30, 46E35, 35D05, 35B50.

\section{Introduction}\label{intro}

A theoretical study of most of the well-known algebras of
generalized functions has pointed out two fundamental structures.
The first one is the algebraic structure of a solid factor ring
$\mathcal{C}$ of generalized numbers. The second one is the
topological structure defined by a family $\mathcal{ P}$ of
seminorms, on a locally convex linear space $E$,
 which is also an algebra. These algebras have been denoted by
$\mathcal{ A}(\mathcal{ C},E,\mathcal{ P})$ and one speaks of
$(\mathcal{ C} ,E,\mathcal{ P})$-algebras of generalized objects.
The definition covers most of the well-known algebras of
generalized functions, as for example, the Colombeau simplified
algebra \cite{C2}, Goursat algebras \cite{M4} and asymptotic
algebras \cite{D1}. On the other hand, special choices for $E$,
$\mathcal{ P}$ and $\mathcal{ C}$ also allow the introduction of
some new algebras. One of them is the so-called Egorov extended
algebra, because of the similarity with the Egorov \cite{E1}
algebra of generalized functions. We have been interested in
working within the framework of the so-called generalized Sobolev
algebras based on the classical Sobolev spaces. As $E$ is a
differential algebra, the main interest of these algebras is to
give a framework which is well suitable to solve many non linear
differential problems with irregular data. The method is based on
the extension of a mapping from $(E_1,\mathcal{ P}_1)$ into
$(E_2,\mathcal{ P}_2)$ to a mapping from $\mathcal{ A}(\mathcal{
C}_1,E_1,\mathcal{ P}_1)$ into $\mathcal{ A} (\mathcal{
C}_2,E_2,\mathcal{ P}_2)$. This method has been introduced, in the
framework of asymptotic algebras, by A. Delcroix and D.
Scarpalezos \cite{D1}, and used, in the framework of $(\mathcal{
C},E,\mathcal{ P})$-algebras, to solve a non linear Dirichlet
problem \cite{M3} and a non linear Neumann problem \cite{M2}, both
with irregular data by  J.-A. Marti and S. P. Nuiro. In this
paper, our goal is to lift up the generalized Sobolev algebras, by
giving more clear definitions of all the statements and general
results in this framework, in order to work more easily with these
algebras. We introduce the first example of ordered generalized
Sobolev algebras, which allows us to pose and eventually solve an
obstacle problem with irregular data. We also point out some
sufficient properties for the existence of an embedding of some
space into a generalized Sobolev algebra.
In the framework of generalized Sobolev algebra, we are able to
solve
 a non linear degenerated Dirichlet problem \cite{M3} with weaker
assumptions.

Consider $\Om$ an open bounded domain of $\R^d$ ($d\in\N^*$) with
a lipschitz continuous boundary $\partial\Om$, we can state this
formal problem :
$$
\left(\mathbf{P}\right)\quad\left\{
\begin{array}{rcl}
-\Delta \Phi\left(u\right)+u=f & in &\Omega,\\
u=g &  on&\partial\Omega,
\end{array}
\right.
$$
where $f$ and $g$ are non smooth functions defined on $\Om$ and
$\partial\Om$ respectively, $\Phi$ an increasing real-valued
differentiable function defined on $\R$ so that $\Phi'$ is a
continuous bounded function that can vanish on a finite set of
discrete points of $\R$. This is a quasilinear diffusion type
problem, with non homogeneous Dirichlet condition on the boundary.
One can remark that  the formal second order differential operator
$ \mathcal{L}=-\dvg\left(\Phi'(.)\grad_{x}\right)+ I_d$ is a
degenerated one, because $\Phi'$ can vanish. Thus,
$\left(\mathbf{P}\right)$ is a Dirichlet nonlinear elliptic
degenerated problem. In order to solve this problem, we introduce
an auxiliary problem by using an artificial viscosity
regularization depending on a parameter $\veps$.

\section{Special types of generalized algebras}\label{stoga}

\subsection{Definitions}
\label{n11111312}

Let us, first, state that $\K$ is $\R$ or $\C$, and
$\un=\left(\un_{\veps}\right)_{\veps}$ where
 $\un_{\veps}=1\hbox{ for all }\veps.$
The generalized algebras constructed from $E$, a normed
$\K$-algebra, are particular case of
$(\mathcal{C},E,\mathcal{P})$-algebras \cite{M1},
\cite{M3}, \cite{M2}, \cite{M4}.\\

Consider a subring $A$ of the ring $\K^{]0,1]}$ so that $\un\in
A$, and which, as a ring, is solid (with compatible lattice
structure) in the following sense :
\begin{definition}
A is said to be {\em solid} if from $(s_{\veps})_{\veps}\in A$ and
$\abs{t_{\veps}}\leq \abs{s_{\veps}}$ for each $\veps\in ]0,1]$ it
follows that $(t_{\veps})_{\veps}\in A$.
\end{definition}
We also consider an ideal $I_A$ of $A$ which is solid as well, and
so that
\begin{equation}
\forall \left(r_{\veps}\right)_{\veps}\in I_A,\quad\lim_{\veps\to
0}r_{\veps}=0.
\end{equation}
Then, we introduce the factor ring $\mathcal{ C}=A/I_A$, which is
called a \textit{ring of generalized numbers}.


\begin{definition}
Let $E$ be a normed algebra. We shall call N-generalized algebra
all factor algebra
$$\mathcal{ A}(\mathcal{ C},E)=\mathcal{ H}_A\left(E\right)/\mathcal{ I}_{I_A}\left(E\right),$$
where $$\mathcal{ H}_A\left(E\right)=\{(u_{\veps})_{\veps}\in
E^{]0,1]}~/~(\Vert u_{\veps}\Vert_E)_{\veps}\in A^+\}$$ and
$$\mathcal{ I}_{I_A}\left(E\right)=\{(u_{\veps})_{\veps}\in
E^{]0,1]}~/~(\Vert u_{\veps}\Vert_E)_{\veps}\in I_A^+\},$$ when
$\Vert\cdot\Vert_E$ is the norm on $E$,
$A^+=\{(r_{\veps})_{\veps}\in A ~/~\forall \veps>0,\;
r_{\veps}\in\R_+\}$ and $I_A^+=\{(r_{\veps})_{\veps}\in I_A
~/~\forall \veps>0,\; r_{\veps}\in\R_+\}$. Its ring of generalized
numbers is defined as the ring
$$\mathcal{ H}_A\left(\K\right)/\mathcal{ I}_{I_A}\left(\K\right)=\mathcal{ C}=A/I_A.$$
\end{definition}

\begin{remark}
 We remark that the notation
is $\mathcal{ A}(\mathcal{ C},E)$ instead of $\mathcal{
A}(\mathcal{ C},E,\mathcal{P})$ since the family $\mathcal{P}$ is
reduced to one single element. The algebra
$\mathcal{A}(\mathcal{C},E)$ is also a vector space on the field
$\K$.
\end{remark}

\begin{example}
\label{n16112055} With
$$
I_A=\left\{r=(r_\veps)_{\veps}\in\R^{]0,1]}\medspace|\medspace\forall
k\in
    \N^*,\quad \abs{r_\veps}=O\left(\veps^k\right)\right\}
    $$
and
$$
A=\left\{r=(r_\veps)_{\veps}\in\R^{]0,1]}\medspace|\medspace\exists
k\in\Z,
    \quad \abs{r_\veps}=O\left(\veps^k\right)\right\},
    $$
we obtain a polynomial growth type N-generalized algebra.
\end{example}

\begin{example}
We take
$$
I_A=\left\{r=(r_\veps)_{\veps}\in\R^{]0,1]}\medspace|\medspace\exists
    \veps_0\in ]0,1],\quad\forall\veps\in]0,\veps_0],\quad
    r_\veps=0\right\},
    $$
and $A=\R^{]0,1]}$. With such $A$ and $I_A$, we obtain another
N-generalized algebra.
\end{example}

\begin{example}
When $E$ is a Sobolev algebra (that is, for example, on the form
$W^{m+1,p}(\Om)\cap W^{m,\infty}(\Om)$, with $m\in ]0,+\infty[$,
$p\in[1,+\infty[$ and $\Om$ an open subset of $\R^d$ ($d\in
\N^*$)), respectively a Banach algebra, we will speak about
generalized Sobolev algebra, respectively generalized Banach
algebra, instead of N-generalized algebra.
\end{example}

\subsection{Embeddings and weak equalities}\label{eawe}

In the following paragraph, we are going to show a way to embed
$E$ into $\mathcal{A}(\mathcal{C},E)$.

\begin{proposition}
\label{N3262258} The mapping $i_0$ defined on $E$, by :
$$\forall u\in E,\quad i_0(u)=cl\left(u\un_\veps\right)_{\veps},
    $$
is linear and one-to-one from $E$ into
$\mathcal{A}(\mathcal{C},E)$.
\end{proposition}

\bproof For every $u\in E$, we have :
$\left(\norm{u\un_\veps}_E\right)_{\veps}=\norm{u}_E\un.
    $
Furthermore, as $\norm{u}_E\in\K$ and $\un\in A$, there exists
$\lambda\in\N$ so that
$$\forall\veps,\quad\norm{u_{\veps}}_E\leq\lambda\un_{\veps},
    $$
and obviously $\lambda\un\in A^+$. As a consequence of the solid
property which implies that $\left(
u_{\veps}\right)_{\veps}\in\mathcal{H}_A\left(E\right)$, we have
$i_0(u)\in\mathcal{A}(\mathcal{C},E)$. It can easily be proved
that $i_0$ is linear and one-to-one. \eproof

\begin{definition}
The mapping $i_0$ from $E$ into $\mathcal{A}(\mathcal{C},E)$,
defined in proposition \ref{N3262258}, will be the so-called
trivial embedding of $E$ into $\mathcal{A}(\mathcal{C},E)$.
\end{definition}

We can also embed some topological vector space into $\mathcal{
A}(\mathcal{C},E)$. Let $(G,\mathcal{T})$ be a Hausdorff
topological vector space so that there exists a continuous linear
mapping $j$ from $\left(E,\norm{.}_E\right)$ into
$\left(G,\mathcal{T}\right)$.

\begin{definition}
$T\in G$ and $U=cl\left(u_{\veps}\right)_{\veps}\in\mathcal{
A}(\mathcal{ C},E)$ are $(G,\mathcal{T})$-associated if
$$j\left(u_{\veps}\right)\to T\hbox{ in }\left(G,
    \mathcal{T}\right)\hbox{ as }\veps\to 0.
    $$
It will be denoted by $U\overset{G,\mathcal{T}}{\sim}T.$
\end{definition}

\begin{remark}
This definition does not depend on the chosen representative of
$U$. Indeed, let $\left(e_{\veps}\right)_{\veps}\in\mathcal{ I}_{
I_A}\left(E\right).$ Therefore, $\dsp\lim_{\veps\to
0}\norm{e_{\veps}}_E=0,$ which means that $e_{\veps}\to 0\hbox{ in
}\left(E,\norm{.}_E\right)\hbox{ as }\veps\to 0.
    $
Consequently, we have $j\left(e_{\veps}\right)\to 0 \hbox{ in
}\left(G,\mathcal{T}\right)$ as $\veps\to 0.$
\end{remark}

\begin{definition}
Assume that $U=cl\left(u_{\veps}\right)_{\veps},
V=cl\left(v_{\veps}\right)_{\veps}\in \mathcal{ A}(\mathcal{
C},E)$. We shall say that $U$ and $V$ are $(G,\mathcal{T})$-weakly
equals if
$$\left(U-V\right)\overset{G,\mathcal{T}}{\sim}0.
    $$
It will be denoted by $U\overset{G,\mathcal{T}}{\simeq}V.$
\end{definition}

\begin{proposition}
Assume that for every $T\in G$, there exists
$\left(u_{\veps}\right)_{\veps}\in \mathcal{H}_A\left(E\right)$,
so that
$$j\left(u_{\veps}\right)\to T\hbox{ in }(G,\mathcal{T}),
    \hbox{ as }\veps\to 0.
$$
Then, there exists, at least, an embedding $i_{G}$ from
$(G,\mathcal{T})$ into the N-generalized algebra $\mathcal{
A}(\mathcal{ C},E)$. Furthermore, if, for all $v\in E$, there
exists
 $\left(u_{\veps}\right)_{\veps}\in\mathcal{H}_A\left(E\right)$,
so that $\left(u_{\veps}-v\right)_{\veps}\in\mathcal{
I}_{I_A}\left(E\right)$, then
\begin{equation}
\label{N06040014} \forall u\in E,\quad \left(i_{G}\circ
j\right)(u)\overset{G,\mathcal{T}}{\simeq}i_0(u).
\end{equation}
\end{proposition}

\bproof For every $T\in G$, there exists
$\left(u_{\veps}\right)_{\veps}\in \mathcal{ H}_A\left(E\right)$,
so that
$$j\left(u_{\veps}\right)\to T\hbox{ in }\left(G,\mathcal{T}\right)\hbox{ as }
    \veps\to 0.
    $$
Let us state $i_{G}(T)=cl\left(u_{\veps}\right)_{\veps}$. The
mapping $i_{G}$ from $G$ into $\mathcal{ A}(\mathcal{ C},E)$ is
 obviously linear. Let us prove that $i_{G}$ is one-to-one. If
$i_{G}(T)=0$ in $\mathcal{ A}(\mathcal{ C},E)$ then
$$i_{G}(T)=cl\left(e_{\veps}\right)_{\veps}\hbox{ for }
    \left(e_{\veps}\right)_{\veps}\in\mathcal{ I}_{I_A}\left(E\right).
    $$
We have $e_{\veps} \rightarrow 0$\hbox{ in
}$\left(E,\norm{.}_E\right)$ which implies that $j\left(e_{\veps}
\right)\to 0\hbox{ in }\left(G,\mathcal{T}\right)$, whenever
$\veps\to 0$. This leads to $T=0$ in $G$, because
$\left(G,\mathcal{T}\right)$ is a Hausdorff space. The second
property is obvious. \eproof

\begin{remark}
If there exists another such embedding $i'_{G}$ from
$(G,\mathcal{T})$ into the N-generalized algebra
$\mathcal{A}(\mathcal{C},E)$ then
$$
\forall T\in G,\quad
i_{G}(T)\overset{G,\mathcal{T}}{\simeq}i'_{G}(T).
$$
\end{remark}

\begin{example}\label{egalfaibl}
Let $j$ be the canonical embedding of
$\left(L^{\infty}(\Om),\norm{.}_{L^{\infty}(\Om)}\right)$ in
$\left(H^{-2}(\Om),\sigma\left(H^{-2}(\Om),H^{2}_0(\Om)\right)\right)$,
where $\sigma\left(H^{-2}(\Om),H^{2}_0(\Om)\right)$ denotes the
weak topology on $H^{-2}(\Om)$. We will say that $T\in
H^{-2}(\Om)$ and $\mathcal{U}=cl(u_{\veps})_{\veps}\in
\mathcal{A}(\mathcal{C},L^{\infty}(\Om))$ are
$H^{-2}(\Om)$-associated if
$$j(u_{\veps})\to T \mbox{ in } \left(H^{-2}(\Om),\sigma\left(H^{-2}(\Om),H^{2}_0(\Om)\right)\right), \mbox{ as }\veps\to 0,$$
and we will denote $\mathcal{U}\overset{2}{\sim} T$. Moreover, we
will say that $\mathcal{U},\mathcal{V}\in
\mathcal{A}(\mathcal{C},L^{\infty}(\Om))$ are $H^{-2}(\Om)$-weakly
equals if $\mathcal{U}-\mathcal{V}\overset{2}{\sim}0$ and we will
denote $\mathcal{U}\overset{2}{\simeq}\mathcal{V}$.
\end{example}

\subsection{Mapping on N-generalized algebra}\label{monga}

The idea of extension of mapping has been introduced by A.
Delcroix and D. Scarpalezos \cite{D1}, in the framework of
asymptotic algebras. But it is, in fact, a particular case of
definition of mapping on $\mathcal{A}(\mathcal{C},E)$-algebras.

If $\theta=\left(\theta_{\veps}\right)_{\veps}$ is a family of
mappings from a normed algebra $(E,\norm{.}_E)$ into a normed
algebra $(F,\norm{.}_F)$, one can view $\theta$ as a mapping from
the N-generalized algebra $\mathcal{A}(\mathcal{C},E)$ into the
N-generalized algebra $\mathcal{A}(\mathcal{D},F)$, where we have
set $\mathcal{C}=A/I_{A}$ and $\mathcal{D}=B/I_{B}$ when $A,I_A,B$
and $I_B$ are as in \S 2.1. One remarks that the extension theorem
of A. Delcroix and D. Scarpalezos \cite{D1} deals with the case
where $\theta=\left(\theta\right)_{\veps}$.

\begin{theorem}
Let $E$ and $F$ be two normed algebras and
$(\theta_{\veps})_{\veps}$ a family of applications of $E$ in $F$.
We assume that
\begin{enumerate}
\item $A\subset B$ and $I_A\subset I_B$, \item there exists a
family of polynomial functions $(\Psi_{\veps})_{\veps}$ of one
variable with coefficients in $A_+$ so that $$\forall
\veps>0~,~\forall x\in E~,~\Vert
\theta_{\veps}(x)\Vert_F\le\Psi_{\veps}(\Vert x\Vert_E),$$ \item
there exists two families of polynomial functions
$(\Psi_{\veps}^1)_{\veps}$ and $(\Psi_{\veps}^2)_{\veps}$ of one
variable with coefficients in $A_+$ so that $\Psi_{\veps}^2(0)=0$
for all $\veps>0$, and
$$\forall \veps>0~,~\forall x,\xi\in
E~,~\Vert\theta_{\veps}(x+\xi)-
\theta_{\veps}(x)\Vert_F\le\Psi_{\veps}^1(\Vert
x\Vert_E)\Psi_{\veps}^2(\Vert \xi\Vert_E).$$
\end{enumerate}
Then there exists an application $\Theta$ :
$\mathcal{A}$$($$\mathcal{C}$$,E)\rightarrow$$\mathcal{A}$$($$\mathcal{D}$$,F)$,
which associates $cl(\theta_{\veps}(x_{\veps}))_{\veps}$ with
$cl(x_{\veps})_{\veps}$.
\end{theorem}

\bproof First, let $(x_{\veps})_{\veps}$ be in
$\mathcal{H}$$_A(E)$ and let us show that
$(\theta_{\veps}(x_{\veps}))_{\veps}$ is in $\mathcal{H}$$_B(F)$.
We have $(\Vert x_{\veps}\Vert_E)_{\veps}$ in $A_+$ so
$(\Psi_{\veps}(\Vert x_{\veps}\Vert_E))_{\veps}$ is also in $A_+$,
since $(\Psi_{\veps})_{\veps}$ has coefficients in $A_+$. Thus
$(\Vert\theta_{\veps}(x_{\veps})\Vert_F)_{\veps})$ belongs to
$A_+\subset B_+$, due to (1) and (2), which implies what we want.
Then, let $(i_{\veps})_{\veps}$ be in $\mathcal{I}$$_{I_A}(E)$ and
let us show that
$(\theta_{\veps}(x_{\veps}+i_{\veps})-\theta_{\veps}(x_{\veps}))_{\veps}$
is in $\mathcal{I}$$_{I_B}(F)$. Since $(\Vert
x_{\veps}\Vert_E)_{\veps}$ and $(\Vert i_{\veps}\Vert_E)_{\veps}$
are respectively in  $A_+$ and $I_A^+$ then $(\Psi_{\veps}^1(\Vert
x_{\veps}\Vert_E))_{\veps}$ and $(\Psi_{\veps}^2(\Vert
i_{\veps}\Vert_E))_{\veps}$ are respectively in $A_+$ and $I_A^+$
, since, for $i\in\{1,2\},$ $(\Psi_{\veps}^i)_{\veps}$ has
coefficients in $A_+$. Then, $(\Psi_{\veps}^1(\Vert
x_{\veps}\Vert_E)\Psi_{\veps}^2(\Vert i_{\veps}\Vert_E))_{\veps}$
is in $I_A^+$. Thus
$(\Vert\theta_{\veps}(x_{\veps}+i_{\veps})-\theta_{\veps}(x_{\veps})\Vert_F)_{\veps}$
belongs to $I_A^+\subset I_B^+$, due to (1) and (3), which implies
the required result. \eproof

As a consequence, we obtain the following result.

\begin{proposition}
Assume that $A\subset B$ and $I_{A}\subset I_{B}$. If
$(\theta_{\veps})_{\veps}$ is a family of continuous linear
mappings from a normed algebra $E$ into a normed algebra $F$, then
$(\theta_{\veps})_{\veps}$ also defines a mapping $\Theta$ from
$\mathcal{A}(\mathcal{C},E)$ into $\mathcal{A}( \mathcal{D},F)$.
\end{proposition}

\begin{example}
Let $\Om$ be an open subset of $\R^d$ and $E=\Hun\cap\Lp{\infty}$
with $\norm{u}_{E}= \norm{u}_{\Lp{\infty}}+\norm{u}_{\Hun}$. The
canonical embedding $i:u\mapsto u$ is continuous as well as linear
from the Banach algebra $E$ into the Banach algebra $\Lp{\infty}$.
Obviously, the mapping $i$ verifies all the assumptions of the
previous proposition; this is why we can define its extension
$\mathcal{I}$ as a mapping from
$\mathcal{A}\left(\mathcal{C},E\right)$ into
$\mathcal{A}\left(\mathcal{C}, \Lp{\infty}\right)$.
\end{example}

In the same way, one can prove that :

\begin{proposition}
Assume that $\left(\theta_{\veps}\right)_{\veps}$ is a family of
mappings from a normed algebra $E$ into the topological field
$(\K,\abs{.})$, so that
\begin{itemize}
\item there exists a family of polynomial functions
$(\Psi_{\veps})_{\veps}$ of one variable with coefficients in
$A_+$ so that $$\forall \veps>0~,~\forall x\in E~,~\abs{
\theta_{\veps}(x)}\le\Psi_{\veps}(\Vert x\Vert_E),$$ \item there
exists two families of polynomial functions
$(\Psi_{\veps}^1)_{\veps}$ and $(\Psi_{\veps}^2)_{\veps}$ of one
variable with coefficients in $A_+$ so that $\Psi_{\veps}^2(0)=0$
for all $\veps>0$, and
$$\forall \veps>0~,~\forall x,\xi\in
E~,~\abs{\theta_{\veps}(x+\xi)-
\theta_{\veps}(x)}\le\Psi_{\veps}^1(\Vert
x\Vert_E)\Psi_{\veps}^2(\Vert \xi\Vert_E).$$
\end{itemize}
Then there exists an application $\Theta$ :
$\mathcal{A}$$($$\mathcal{C}$$,E)\rightarrow$ $\mathcal{C}$, which
associates $cl(\theta_{\veps}(x_{\veps}))_{\veps}$ with
$cl(x_{\veps})_{\veps}$.
\end{proposition}

\begin{remark}
If $\theta$
is a continuous linear mapping from a normed algebra
$(E,\norm{.}_E)$ into the topological field $(\K,\abs{.})$, then
$\theta$ also defines a mapping, denoted by $\Theta$, from
$\mathcal{A}(\mathcal{C},E)$ into the factor ring
$\mathcal{C}=A/I_A$.
\end{remark}

\subsection{An example of ordered generalized Sobolev
algebra}\label{aeoogsa}

Consider $A$ and $I_A$ as in \S \ref{n11111312}, the Sobolev
algebra $\Lp{\infty}$, endowed with its usual topology, with $\Om$
an open bounded subset of
    $\R^d$.
Thus, we can consider the algebra
$\mathcal{A}(\mathcal{C},\Lp{\infty})$. It is easy to prove, by
means of theorem 1, that the mapping
$$
\begin{array}{rcl}
p:  \Lp{\infty} & \rightarrow & \Lp{\infty} \\
     u          & \mapsto     & u^+=\sup{\{u,0\}}=\frac 12 (u+\abs{u})
\end{array}
$$
can be extended as a mapping $\mathcal{P}$ from
$\mathcal{A}(\mathcal{C},\Lp{\infty})$ into itself, defined by :
$$
\forall
U=cl(u_{\veps})_{\veps}\in\mathcal{A}(\mathcal{C},\Lp{\infty}),\quad
    \mathcal{P}(U)=cl(p(u_{\veps}))_{\veps},
$$
due to the following relation :
$$
\forall r,s\in\R,\quad\abs{(r+s)^+-r^+}\leq\abs{s}.
$$

We are now able to state the following result:

\begin{proposition}
The generalized Sobolev algebra
$\mathcal{A}(\mathcal{C},\Lp{\infty})$ is partially ordered by the
following binary relation :
$$
\forall U,V\in\mathcal{A}(\mathcal{C},\Lp{\infty}),\quad U\leq V
\Longleftrightarrow \mathcal{P}\left(U-V\right)=0.
$$
\end{proposition}

\bproof Obviously, the relation $\leq$ is reflexive, then we have
to prove, for $U,V,W\in\mathcal{A} (\mathcal{C},\Lp{\infty})$,
that :
\begin{eqnarray}
U\leq V\textit{ and }V\leq U &\Rightarrow& U=V;\label{n11111732}\\
U\leq V\textit{ and }V\leq W &\Rightarrow& U\leq
W.\label{n11111745}
\end{eqnarray}
We state $U=cl(u_{\veps})_{\veps}$, $V=cl(v_{\veps})_{\veps}$ and
$W=cl(w_{\veps})_{\veps}$.
\newline
\textit{Proof of property (\ref{n11111732}) :} If $U\leq V\text{
and }V\leq U$ then, there exists
$\left(\varphi_{\veps}\right)_{\veps}$ and
$\left(\psi_{\veps}\right)_{\veps}$ in
$\mathcal{I}_{I_A}\left(\Lp{\infty}\right)$ so that
$\left(u_{\veps}-v_{\veps}\right)^{+}=\varphi_{\veps}$ and
$\left(v_{\veps}- u_{\veps}\right)^{+}=\psi_{\veps}$. As,
$$
u_{\veps}-v_{\veps}=\left(u_{\veps}-v_{\veps}\right)^{+}-\left(v_{\veps}-u_{\veps}\right)^{+}=
\varphi_{\veps}-\psi_{\veps},
$$
it follows that $ \left(u_{\veps}-v_{\veps}\right)_{\veps}
=\left(\varphi_{\veps}-\psi_{\veps}\right)_{\veps}\in
\mathcal{I}_{I_A}\left(\Lp{\infty}\right), $ whence $U=V$.\newline
\textit{Proof of property (\ref{n11111745}) :} If $U\leq V\text{
and }V\leq W$ then we have
$$
\left(\norm{\left(u_{\veps}-v_{\veps}\right)^{+}}_{\Lp{\infty}}\right)_\veps\in
I_A\,\,\,,
\left(\norm{\left(v_{\veps}-w_{\veps}\right)^{+}}_{\Lp{\infty}}\right)_\veps\in
I_A.
$$
By means of the solid property, we deduce, from the following
inequality :
$$
\norm{\left(u_{\veps}-w_{\veps}\right)^{+}}_{\Lp{\infty}}\leq\norm{\left(u_{\veps}-v_{\veps}
\right)^{+}}_{\Lp{\infty}}+\norm{\left(v_{\veps}-w_{\veps}\right)^{+}}_{\Lp{\infty}},
$$
that
$\left(\left(u_{\veps}-w_{\veps}\right)^+\right)_{\veps}\in\mathcal{I}_{I_A}\left(\Lp{\infty}
\right), $ which yields $\mathcal{P}(U-W)=0$, that is to say
$U\leq W$. \eproof

\begin{proposition}
For all $u,v\in\Lp{\infty}$, we have $i_0(u)\leq i_0(v)$ if, and
only if, $u\leq v$ in $\Lp{\infty}$, that is $u\leq v$ almost
everywhere in $\Om$.
\end{proposition}

\bproof If $i_0(u)\leq i_0(v)$ then
$\mathcal{P}\left(i_0(u)-i_0(v)\right)=0$. Consequently, there
exists
$\left(\varphi_{\veps}\right)_{\veps},(e_{\veps})_{\veps}\in\mathcal{I}_{I_A}\left(\Lp{\infty}
\right)$ so that $\left(u-v+e_{\veps}\right)^{+}=\varphi_{\veps}$,
since we have $u_{\veps}-v_{\veps}=u-v+e_{\veps}$ for all $\veps$.
Taking into account that
$$
\varphi_{\veps}\to 0\quad and \quad e_{\veps}\to
0\quad\textit{in}\,\,\Lp{\infty}\,\, \textit{, as}\,\,\veps\to 0,
$$
it may be seen that $
\left(u-v+e_{\veps}\right)^{+}=\varphi_{\veps}\to
0\quad\textit{a.e. in }\Om, $ whence
$\left(u-v\right)^{+}=0\quad\text{a.e. in }\Om$, since one can
easily prove that
$$
\left(u-v+e_{\veps}\right)^{+}\to(u-v)^{+}\quad\textit{in}\,\,\Lp{\infty}\,\,
\textit{, as}\,\,\veps\to 0.
$$
It means that $u\leq v\quad\text{a.e. in }\Om.$

Conversely, if $u\leq v\;\text{a.e. in }\Om$ then
$(u-v)^+=0\;\text{a.e. in }\Om$. By definition of $\mathcal{P}$
and $i_0$, this leads to $i_0(u)\leq i_0(v)$. \eproof

\begin{proposition}
Let $u\in \Lp{\infty}$ and
$U\in\mathcal{A}(\mathcal{C},\Lp{\infty})$. If
$U\overset{\Lp{1}}{\sim}u$ (here $L^1(\Om)$ is endowed with its
usual topology) and $U\leq 0$ then $u\leq 0$ a.e. in $\Om$.
\end{proposition}

\bproof We set $U=cl(u_{\veps})_{\veps}$. Since
$U\overset{\Lp{1}}{\sim}u$ then, as $\veps$ goes to $0$,
$u_{\veps}\to u\text{ in }\Lp{1}$, which gives $u_{\veps}^+\to
u^+\text{ in } \Lp{1}$, by means of the Lebesgue dominated
convergence theorem. Since $U\leq 0$ then
$\mathcal{P}\left(U\right)=0$. Consequently, there exists a
sequence of functions
$\left(\varphi_{\veps}\right)_{\veps}\in\mathcal{I}_{I_A}\left(\Lp{\infty}\right)$
so that $u_{\veps}^{+}=\varphi_{\veps}$ for all $\veps$. Taking
into account that
$$
\varphi_{\veps}\to 0\quad\textit{in}\,\,\Lp{\infty}\,\,\textit{,
as}\,\,\veps\to 0,
$$
we find that $ u_{\veps}^{+}=\varphi_{\veps}\to 0\quad\textit{a.e.
in }\Om, $ whence $u^{+}=0$ a.e. in $\Om,$ which implies $u\leq 0$
a.e. in $\Om.$ \eproof

\section{Solution of the nonlinear degenerate Dirichlet
problem}\label{sotnddp}

After having solved the auxiliary problem by using an artificial
viscosity regularization depending on a parameter $\veps$, we
solve our main problem $\left(\mathbf{P}\right)$ (see section
\ref{intro}), in a generalized Sobolev algebra with the classical
equality and with the weak one defined in example \ref{egalfaibl}.
Then we perform a little qualitative study of the solution.

\subsection{The regularized Dirichlet problem}\label{solreg}

Let us set $$\mathcal{V}\displaystyle
_A^+=\left\{(r_{\veps})_{\veps}\in
A^+~/~\forall\veps>0,\;r_{\veps}\in\ ]0,1];~\lim_{\veps\to
0}r_{\veps}=0~;~\left(\frac{1}{r_{\veps}}\right)_{\veps}\in
A^+\right\}.$$ Assume that $\mathcal{V}$$\displaystyle
_A^+\not=\emptyset$ and then, for all $(r_{\veps})_{\veps}$ in
$\mathcal{V}$$_A^+$, set $\Phi_{\veps}=\Phi+r_{\veps}Id$. This
section consists in proving the following proposition :
\begin{proposition}
If $f\in L^{\infty}(\Om)$ and $g\in L^{\infty}(\partial\Om)$ then
there exists one, and only one, function $u\in H^1(\Om)\cap
L^{\infty}(\Om)$ solution of the regularized problem
$$
\left(\mathbf{P_{\veps}}\right)\quad\left\{
\begin{array}{rcl}
-\Delta \Phi_{\veps}\left(u\right)+u=f & in &\Omega,\\
u=g &  on&\partial\Omega.
\end{array}
\right.
$$
\end{proposition}

\bproof
This proof goes in three steps.\\

1)Maximum's principle\\  We are going to prove that if $u\in
H^1(\Omega)$ is a solution of this problem then
$$m\le u\le M\mbox{ a.e. in } \Omega,$$
with $\displaystyle
m=\min\{\inf_{\Omega}f~,~\inf_{\partial\Omega}g\}\mbox{ and
}M=\max\{\sup_{\Omega}f~,~\sup_{\partial\Omega}g\},$ which means
that
$u$ belongs to $L^{\infty}(\Omega)$.\\
Indeed, for such a $u$, we have, for all $v$ in $H^1_0(\Omega)$
$$\int_{\Omega}\grad\Phi_{\veps}(u)\grad v dx + \int_{\Omega}u v dx =
\int_{\Omega} fv dx,$$ where $dx$ denotes the Lebesgue measure on
$\Omega$. Let us consider the function $\displaystyle
v=(\Phi_{\veps}(u)-\Phi_{\veps}(M))^+$ then $v$ is in
$H^1_0(\Omega)$, so
\begin{multline*}
\int_{\Omega}(\grad(\Phi_{\veps}(u)-\Phi_{\veps}(M))^+)^2 dx + \int_{\Omega}u (\Phi_{\veps}(u)-\Phi_{\veps}(M))^+ dx =\\
\int_{\Omega} f(\Phi_{\veps}(u)-\Phi_{\veps}(M))^+ dx,
\end{multline*} since
$\Phi_{\veps}(M)$ is a constant. Consequently,
\begin{multline*}
\Vert (\Phi_{\veps}(u)-\Phi_{\veps}(M))^+\Vert^2_{H^1_0(\Omega)}
=\int_{\Omega} (f-M)(\Phi_{\veps}(u)-\Phi_{\veps}(M))^+ dx \\-
\int_{\Omega} (u-M)(\Phi_{\veps}(u)-\Phi_{\veps}(M))^+ dx.
\end{multline*}
By definition of $M$, the first integral is negative and, since
the functions $Id$ and $\Phi_{\veps}$ are increasing, the second
one is non negative. Then
$$\Vert
(\Phi_{\veps}(u)-\Phi_{\veps}(M))^+\Vert^2_{H^1_0(\Omega)}\le 0,$$
that is $\Phi_{\veps}(u)\le\Phi_{\veps}(M)$ a.e. in $\Omega$,
which implies the first part of the required result, since
$\Phi_{\veps}$ is an increasing function. For the second part, we
use a similar
method by taking $\displaystyle v=(\Phi_{\veps}(u)-\Phi_{\veps}(m))^-$.\\

2)Existence of a solution in $H^1(\Omega)$\\
This result is obtained by using the Schauder's fixed point
theorem related to a weakly sequentially continuous mapping from a
reflexive and separable Banach space into itself. Let us consider
$w_0\in H^1(\Omega)$ the unique solution of the following linear
Dirichlet problem :
$$\left\{
\begin{array}{rll}
-\Delta w_0&=&0 \mbox{ in } \Omega,\\
w_0&=&g \mbox{ on } \partial\Omega.
\end{array}
\right.$$ Then a solution of the regularized problem is of the
form $w_0+w$, with $w\in H^1_0(\Omega)$ and for all $v$ in
$H^1_0(\Omega)$, one has
$$\int_{\Omega}\Phi_{\veps}'(w_0+w)\grad w_0 \grad v dx + \int_{\Omega}\Phi_{\veps}'(w_0+w)\grad w \grad v dx +
     \int_{\Omega}(w_0+w) v dx = \int_{\Omega} fv dx.$$
Consequently, for all $h\in H^1_0(\Omega)$, let us look for $w_h$
in $H^1_0(\Omega)$ so that, for all $v$ in $H^1_0(\Omega)$,
\begin{multline*}
\int_{\Omega}\left\{\tilde{\Phi}_{\veps}'(w_0+h)\grad w_h \grad v
+w_hv\right\}dx\\ = \int_{\Omega}\left\{ fv-w_0v-
\tilde{\Phi}_{\veps}'(w_0+h)\grad w_0 \grad v\right\}dx,
\end{multline*}
where $\tilde{\Phi}_{\veps}$ is defined by
$$\tilde{\Phi}_{\veps}(x)=
\left\{
\begin{array}{ll}
\Phi(m)+r_\veps x& \mbox{ if }x\le m\\
\Phi(x)+r_\veps x& \mbox{ if }x\in]m,M[\\
\Phi(M)+r_\veps x& \mbox{ if }x\ge M.
\end{array}
\right.$$ The existence and uniqueness of $w_0$ and $w_h$ are
ensured by the Lax-Milgram's theorem. Moreover, for the
test-function $v=w_h$, we get
\begin{multline*}
\int_{\Omega}\left\{\tilde{\Phi}_{\veps}'(w_0+h)\vert\grad
w_h\vert^2+\abs{w_h}^2\right\} dx= \int_{\Omega} f w_h dx \\-
     \int_{\Omega}\left\{\tilde{\Phi}_{\veps}'(w_0+h)\grad w_0 \grad w_h - w_0w_h\right\} dx.
\end{multline*}
Meanwhile,
$$\int_{\Omega}\left\{\tilde{\Phi}_{\veps}'(w_0+h)\vert\grad w_h\vert^2+\abs{w_h}^2\right\}
                                  dx\ge r_{\veps}\int_{\Omega}\left\{\vert\grad w_h\vert^2 +\abs{w_h}^2\right\}dx,$$
$$\int_{\Omega} f w_h dx\le C(\Omega) \Vert
                     f\Vert_{L^{\infty}(\Omega)}\Vert
                     w_h\Vert_{H^1_0(\Omega)}$$
and
\begin{multline*}
-\int_{\Omega}\left\{\tilde{\Phi}_{\veps}'(w_0+h)\grad w_0 \grad
w_h+w_0w_h\right\} dx \\ \le \dsp C(\Om)(1+r_{\veps}+\Vert
\Phi'\Vert_{L^{\infty}(\R)})\Vert
w_0\Vert_{H^1(\Omega)}\Vert w_h\Vert_{H^1_0(\Omega)}\\
\le\dsp C(\Om)(2+\Vert \Phi'\Vert_{L^{\infty}(\R)})\Vert
w_0\Vert_{H^1(\Omega)}\Vert w_h\Vert_{H^1_0(\Omega)},
\end{multline*}
where $C(\Omega)$ denotes a constant depending on $\Omega$. Thus,
$$\displaystyle \Vert w_h\Vert_{H^1_0(\Omega)}\le\frac{1}{r_{\veps}}C(\Omega)\left[ \Vert
                     f\Vert_{L^{\infty}(\Omega)}+(2+\Vert
\Phi'\Vert_{L^{\infty}(\R)})\Vert w_0\Vert_{H^1(\Omega)}\right]
                     .$$
Noticing that $\Vert w_0\Vert_{H^1(\Omega)}$ depends only on $g$
and $\Omega$ and not on $\veps$, we obtain that $\displaystyle
\Vert w_h\Vert_{H^1_0(\Omega)}\le\frac{C(\Omega,f,g)}{r_{\veps}}$,
which implies that the closed ball $B(0,R_{\veps})$ of center $0$
and radius $R_{\veps}=\frac{C(\Omega,f,g)}{r_{\veps}}$ of the
separable Hilbert space $H^1_0(\Omega)$ is stable by the
application
$$\Pi~:~\begin{array}[t]{ccc}
H^1_0(\Omega)&\rightarrow &H^1_0(\Omega)\\
h            &\mapsto     &w_h.
\end{array}$$
Now we have to prove that for all sequence $(h_n)_n$ of
$B(0,R_{\veps})$ converging weakly to $h$, when $n$ tends to
$+\infty$, the sequence $(\Pi(h_n))_n$ converges weakly to
$\Pi(h)$. Let us consider such a sequence $(h_n)_n$. Since
$(\Pi(h_n))_n$ is bounded, we can extract a subsequence, still
denoted by $(\Pi(h_n))_n$, so that
$$\Pi(h_n)\rightharpoonup \chi \mbox{ in }H^1_0(\Omega).$$
As the imbedding  of $H^1_0(\Omega)$ into $L^2(\Omega)$ is
compact, after another extraction, we have
$$\left\{
\begin{array}{rcl}
\Pi(h_n)&\to& \chi \mbox{ in }L^2(\Omega),\\
h_n&\to &h \mbox{ in }L^2(\Omega) \mbox{ and a.e. in }\Omega.
\end{array}
\right.$$ Since $\tilde{\Phi}_{\veps}'$ is a bounded and piecewise
continuous function and, using the Lebesgue dominated convergence
theorem, we have also
$$\tilde{\Phi}_{\veps}'(w_0+h_n)\to
\tilde{\Phi}_{\veps}'(w_0+h)\mbox{ in }L^2(\Omega).$$ Moreover,
for all $n$ in $\N$ and all $v$ in $H^1_0(\Omega)$, we have
\begin{multline*}
\int_{\Omega}\left\{\tilde{\Phi}_{\veps}'(w_0+h_n)\grad w_{h_n}
\grad v +w_{h_n}v\right\}dx\\ = \int_{\Omega}\left\{ fv-w_0v-
\tilde{\Phi}_{\veps}'(w_0+h_n)\grad w_0 \grad v\right\}dx.
\end{multline*}
Passing to the limit, as $n$ tends to the infinity, in this
previous equality, we obtain that, for all $v$ in $H^1_0(\Omega)$,
\begin{multline*}
\int_{\Omega}\left\{\tilde{\Phi}_{\veps}'(w_0+h)\grad \chi \grad v
+\chi v\right\}dx\\ = \int_{\Omega}\left\{ fv-w_0v-
\tilde{\Phi}_{\veps}'(w_0+h)\grad w_0 \grad v\right\}dx.
\end{multline*}
Meanwhile, for all $h$ in $H^1_0(\Omega)$, there is one and only
one $w_h=\Pi(h)$, so $\Pi(h)=\chi$ and the whole sequence
$(\Pi(h_n))_n$ converges weakly to $\Pi(h)$ in $H^1_0(\Omega)$. We
can now apply the fixed point theorem and conclude that there is
$w$ in $H^1_0(\Omega)$ so that $\Pi(w)=w$. Setting $u=w_0+w$, we
have $u$ in $H^1(\Omega)$ and, for all $v$ in $H^1_0(\Omega)$,
$$\int_{\Omega}\tilde{\Phi}_{\veps}'(u)\grad u \grad v dx +\int_{\Omega}u v dx= \int_{\Omega} fv dx,$$
that is to say that $u$ is solution of
$$\left\{
\begin{array}{rll}
-\Delta \tilde{\Phi}_{\veps}(u)+u&=&f \mbox{ in } \Omega,\\
u&=&g \mbox{ on } \partial\Omega.
\end{array}
\right.$$ Using a method similar to the first step, for this
problem, we can prove that
$$m\le u\le M\mbox{ a.e. in } \Omega,$$
which shows, in fact, that $u$ is solution of the regularized
problem and $u$ belongs to $H^1(\Omega)\cap
L^{\infty}(\Omega)$.\\
Moreover, $$\displaystyle \Vert u\Vert_{H^1(\Omega)}\le
                      \Vert w_0\Vert_{H^1(\Omega)}+\Vert w\Vert_{H^1_0(\Omega)}.$$
But, by definition of $w_0$, we have $\Vert
w_0\Vert_{H^1(\Omega)}\le C(\Omega)\Vert
g\Vert_{L^{\infty}(\partial\Omega)}$ and we prove that
$$\displaystyle \Vert w\Vert_{H^1_0(\Omega)}\le
\frac{1}{r_{\veps}}C(\Omega)\left[\Vert
f\Vert_{L^{\infty}(\Omega)}
      +(2+\Vert
\Phi'\Vert_{L^{\infty}(\R)})\Vert w_0\Vert_{H^1(\Omega)}\right],$$
 so
 \begin{equation}\label{estim}
 \displaystyle \Vert u\Vert_{H^1(\Omega)}\le
      \frac{C(\Omega)}{r_{\veps}}\left[\Vert f\Vert_{L^{\infty}(\Omega)}
      +(2+\Vert
\Phi'\Vert_{L^{\infty}(\R)})\Vert
g\Vert_{L^{\infty}(\partial\Omega)}\right].
 \end{equation}

3)Uniqueness of the solution in $H^1(\Omega)$\\
Let $u_1$ and $u_2$ in $H^1(\Omega)$ be two solutions of the
regularized problem, then for all $v$ belonging to
$H^1_0(\Omega)$, one has
$$\int_{\Omega}\grad(\Phi_{\veps}(u_1)-\Phi_{\veps}(u_2)) \grad v dx + \int_{\Omega}(u_1-u_2) v dx =0.$$
Taking $v=\Phi_{\veps}(u_1)-\Phi_{\veps}(u_2)$, we can write that
$$\Vert \Phi_{\veps}(u_1)-\Phi_{\veps}(u_2)\Vert^2_{H^1_0(\Omega)}+
\int_{\Omega}(u_1-u_2) (\Phi_{\veps}(u_1)-\Phi_{\veps}(u_2)) dx
=0.$$ But it is the sum of two non negative terms, so both are
equal to zero. In particular, $\Vert
\Phi_{\veps}(u_1)-\Phi_{\veps}(u_2)\Vert^2_{H^1_0(\Omega)}=0$,
that is $u_1=u_2$, since $\Phi_{\veps}$ is an injective
function.\eproof

\subsection{Strong solution of the generalized Dirichlet
problem}\label{ssotgdp}

We are going to apply theorem 1 with $E=L^{\infty}(\Omega)\times
L^{\infty}(\partial\Omega)$, $F=H^1(\Omega)\cap
L^{\infty}(\Omega)$ and $$\theta_{\veps}~:~
\begin{array}[t]{ccc}
E     & \rightarrow & F\\
(f,g) & \mapsto     & \theta_{\veps}(f,g)=u_{\veps},
\end{array}$$
where $u_{\veps}$ is the solution of problem
$(\mathbf{P_{\veps}})$. Before, we are going to show the two
following lemmas.

\begin{lemma}
For all $(f,g)$ in $E$ and $u_{\veps}=\theta_{\veps}(f,g)$ in $F$,
we have
$$\displaystyle \Vert u_{\veps}\Vert_F\le
      \frac{C(\Omega,\Vert
\Phi'\Vert_{L^{\infty}(\R)})}{r_{\veps}}\Vert (f,g)\Vert_E.$$
\end{lemma}

\bproof This result is an immediate consequence of inequality
(\ref{estim}) since $\max\{\vert m\vert~,~\vert M\vert\}$ is less
than $\Vert (f,g)\Vert_E=\Vert f\Vert_{L^{\infty}(\Omega)}+\Vert
g\Vert_{L^{\infty}(\partial\Omega)}$. \eproof

\begin{lemma}
For all $(f,g)$, $(\delta,\eta)$ in $E$,
$u_{\veps}=\theta_{\veps}(f,g)$ in $F$ and
$u_{\veps}+\nu_{\veps}=\theta_{\veps}(f+\delta,g+\eta)$ in $F$, we
have
$$\displaystyle \Vert \nu_{\veps}\Vert_F\le
      \frac{C(\Omega,\Vert
\Phi'\Vert_{L^{\infty}(\R)})}{r_{\veps}}\Vert (\delta,\eta)
\Vert_E.$$
\end{lemma}

\bproof By definition of $\theta_{\veps}$, we have
$$\left\{
\begin{array}{rll}
-\Delta \Phi_{\veps}(u_{\veps})+u_{\veps}&=&f \mbox{ in } \Omega,\\
u_{\veps}&=&g \mbox{ on } \partial\Omega,
\end{array}
\right.$$ and $$\left\{
\begin{array}{rll}
-\Delta \Phi_{\veps}(u_{\veps}+\nu_{\veps})+u_{\veps}+\nu_{\veps}&=&f+\delta \mbox{ in } \Omega,\\
u_{\veps}+\nu_{\veps}&=&g +\eta\mbox{ on } \partial\Omega,
\end{array}
\right.$$ so $$\left\{
\begin{array}{rll}
-\Delta \chi_{\veps}(\nu_{\veps})+\nu_{\veps}&=&\delta \mbox{ in } \Omega,\\
\nu_{\veps}&=&\eta \mbox{ on } \partial\Omega,
\end{array}
\right.$$ with
$\chi_{\veps}=\Phi_{\veps}(u_{\veps}+\cdot)-\Phi_{\veps}(u_{\veps})$
which satisfies the same hypothesis as $\Phi_{\veps}$ of section
\ref{solreg}. Consequently, $\nu_{\veps}$ is the solution of a
similar problem as $(\mathbf{P_{\veps}})$ and satisfies an
inequality of the same type as (\ref{estim}), that is
$$\displaystyle \Vert \nu_{\veps}\Vert_F\le
      \frac{C(\Omega,\Vert
\chi'\Vert_{L^{\infty}(\R)})}{r_{\veps}}\Vert (\delta,\eta)
\Vert_E,$$ where $\chi=\chi_{\veps}-r_{\veps} Id$. And inequality
$\Vert \chi'\Vert_{L^{\infty}(\R)}\le\Vert
\Phi'\Vert_{L^{\infty}(\R)}$ implies the required result.\eproof

\begin{theorem}
\label{n0911031130} If $($$\mathcal{F}$$,$$\mathcal{G}$$)$ belongs
to $\mathcal{A}$$($$\mathcal{C}$$,L^{\infty}(\Omega)\times
L^{\infty}(\partial\Omega))$ then there is one, and only one,
generalized function $\mathcal{U}$$=cl(u_{\veps})_{\veps}$,
belonging to $\mathcal{A}$$($$\mathcal{C}$$,H^1(\Omega)\cap
L^{\infty}(\Omega))$, so that
\begin{equation}\label{pbgenreg}
\left\{
\begin{array}{rll}
cl[-\Delta \Phi_{\veps}(u_{\veps})]_{\veps}+\mathcal{{U}}&=&\mathcal{F} \mbox{ in } \mathcal{A}(\mathcal{C},L^{\infty}(\Om)),\\
\Gamma(\mathcal{U})&=&\mathcal{G} \mbox{ in }
\mathcal{A}(\mathcal{C},L^{\infty}(\partial\Om)),
\end{array}
\right.
\end{equation}
where, by definition,
$\Gamma(\mathcal{U})=cl(u_{\veps_{\mid\partial\Om}})_{\veps}=cl(g_{\veps})_{\veps}$,
when ${\mathcal{G}}=cl(g_{\veps})_{\veps}$.
\end{theorem}

\bproof We are going to apply theorem 1 with
$E=L^{\infty}(\Omega)\times L^{\infty}(\partial\Omega)$,
$F=H^1(\Omega)\cap L^{\infty}(\Omega)$ and $$\theta_{\veps}~:~
\begin{array}[t]{ccc}
E     & \rightarrow & F\\
(f,g) & \mapsto     & \theta_{\veps}(f,g)=u_{\veps},
\end{array}$$
where $u_{\veps}$ is the solution of problem
$(\mathbf{P_{\veps}})$. In order to obtain the required result, it
suffices to use the two previous lemmas and apply theorem 1 with
$$\Psi_{\veps}(x)=\frac{C(\Omega,\Vert
\Phi'\Vert_{L^{\infty}(\R)})}{r_{\veps}} x=\Psi_{\veps}^2(x)$$ and
$\Psi_{\veps}^1(x)=1$, for all $x$ in $\R$. The fact that $\Phi'$
is bounded, ensures that $\left(\frac{C(\Omega,\Vert
\Phi'\Vert_{L^{\infty}(\R)})}{r_{\veps}}\right)_\veps$ is in
$A^+$. We set then
$\mathcal{U}$$=\Theta($$\mathcal{F}$$,$$\mathcal{G}$$)=cl(u_{\veps})_{\veps}=cl(\theta_{\veps}(f_{\veps},g_{\veps}))_{\veps}$
when $\mathcal{F}$$=cl(f_{\veps})_{\veps}$ and
$\mathcal{G}$$=cl(g_{\veps})_{\veps}$.\eproof

\subsection{Weak solution of the generalized Dirichlet
problem}\label{wsotgdp}

In this section, we define the notion of weak solution by using
the weak equality defined in example \ref{egalfaibl}.

\begin{theorem}
\label{n0911031140} With the assumptions of theorem
\ref{n0911031130}, if $\mathcal{F}=cl(f_{\veps})_{\veps}$ and
$\mathcal{G}=cl(g_{\veps})_{\veps}$ are such that
\begin{equation}\label{hyp}
\displaystyle \exists (r_{\veps})_{\veps}\in\mathcal{V}_A^+,
\lim_{\veps\to0^+}r_\veps\max\left\{\Vert
g_{\veps}\Vert_{L^{\infty}(\partial\Om)},\Vert
f_{\veps}\Vert_{L^{\infty}(\Om)}\right\}=0,
\end{equation}
then there is one, and only one, generalized function
$\mathcal{U}$ which belongs to
$\mathcal{A}(\mathcal{C},H^1(\Omega)\cap L^{\infty}(\Omega))$ and
such that
\begin{equation}\label{pbgenass}
\left\{
\begin{array}{rcl}
-\Delta \Phi(\mathcal{U})+\mathcal{U}&\overset{2}{\simeq}&\mathcal{F} \mbox{ in } \mathcal{A}(\mathcal{C},L^{\infty}(\Om)),\\
\Gamma(\mathcal{U})&=&\mathcal{G} \mbox{ in }
\mathcal{A}(\mathcal{C},L^{\infty}(\partial\Om)),
\end{array}
\right.
\end{equation}
with $\Delta \Phi(\mathcal{U})=cl(\Delta
\Phi(u_{\veps}))_{\veps}=cl(u_\veps-r_\veps\Delta
u_\veps-f_\veps)_\veps$.
\end{theorem}

\bproof Since $\displaystyle cl[-\Delta
\Phi_{\veps}(u_{\veps})]_{\veps}+\mathcal{{U}}=\mathcal{F} \mbox{
in } \mathcal{A}(\mathcal{C},L^{\infty}(\Om))$, so
$H^{-2}(\Om)$-weakly equal and $\Phi_{\veps}=\Phi+r_{\veps} Id$,
it is sufficient to prove that
$$\displaystyle
Cl(-r_{\veps}\Delta u_{\veps})_{\veps}\overset{2}{\simeq}0.$$ Let
$\varphi$ be in $H^2_0(\Om)$, using Green's formula, one has
$$\displaystyle \int_{\Om}-r_{\veps}\Delta u_{\veps}\varphi dx
\begin{array}[t]{l}
\displaystyle =r_{\veps}\left(\int_{\Om}\nabla
u_{\veps}\nabla\varphi
        dx-\int_{\partial\Om}g_{\veps}\varphi d\nu\right)\\
\displaystyle =r_{\veps}\left(\int_{\partial\Om}u_{\veps}\frac{\partial\varphi}{\partial\nu}d\nu-\int_{\Om}u_{\veps}\Delta\varphi dx -\int_{\partial\Om}g_{\veps}\varphi d\nu\right)\\
\displaystyle =-r_{\veps}\int_{\Om}u_{\veps}\Delta\varphi dx.
\end{array}$$
Consequently, using Cauchy-Schwartz's inequality, one has
$$\left\vert\int_{\Om}-r_{\veps}\Delta u_{\veps}\varphi dx\right\vert
\le r_{\veps}\max(\Vert
    g_{\veps}\Vert_{L^{\infty}(\partial\Om)},\Vert
    f_{\veps}\Vert_{L^{\infty}(\Om)})C(\Om)\norm{\Delta\varphi}_{\Lp{2}}.$$
The assumption (\ref{hyp}) implies that
$$\displaystyle\lim_{\veps\to 0}\int_{\Om}-r_{\veps}\Delta
u_{\veps}\varphi dx=0.$$\eproof
\begin{remark}
This theorem leads us to notice that we can have a Dirac
generalized function in the second member of the problem. Indeed,
a representative of a Dirac generalized function can be : $\forall
x\in\R^d,\;\delta_{\veps}(x)=\veps^{-d}\varphi\left(\veps^{-1}x\right)$,
where $\varphi$ is a compactly supported function defined on
$\R^d$. The hypothesis (\ref{hyp}) is satisfied with: $\forall
q\in\N,\;r_{\veps}=\veps^{d+q}$, and, for example,  we take $A$
and $I_A$ as in example \ref{n16112055}.
\end{remark}

\subsection{Non positive solutions}\label{nps}

In this section, we prove that the solution is non positive, in a
sense to be defined, when the data is. We start by defining what
non positive means here. \bdf An element
$\mathcal{U}\in\mathcal{A}(\mathcal{C},H^1(\Omega)\cap
L^{\infty}(\Omega))$ is said to be non positive if and only if the
corresponding element  $\mathcal{I}\left(\mathcal{U}\right)$, of
the generalized Sobolev algebra
$\mathcal{A}(\mathcal{C},L^{\infty}(\Omega))$, is non positive.
\edf In this definition, $\mathcal{I}$ denotes the extension of
the canonical embedding of $H^1(\Omega)\cap L^{\infty}(\Omega))$
into $L^{\infty}(\Omega)$, introduced in example 5. This mapping
is an embedding of $\mathcal{A}(\mathcal{C},H^1(\Omega)\cap
L^{\infty}(\Omega))$ into
$\mathcal{A}(\mathcal{C},L^{\infty}(\Omega))$.

\bpp With the assumptions of theorem \ref{n0911031140}, if the
generalized functions
$\mathcal{F}$$=cl(f_{\veps})_{\veps}\in\mathcal{A}(\mathcal{C},L^{\infty}(\Omega))$
and
$\mathcal{G}$$=cl(g_{\veps})_{\veps}\in\mathcal{A}(\mathcal{C},L^{\infty}(\partial\Omega))$
are non positive, then
$\mathcal{U}=\Theta(\mathcal{F},\mathcal{G})=cl(\theta_{\veps}(f_{\veps},g_{\veps}))_{\veps}\in\mathcal{A}
(\mathcal{C},H^1(\Omega)\cap L^{\infty}(\Omega))$, the solution to
our main problem, is non positive. \epp

\bproof
 Using the hypothesis on
$\mathcal{F}$, $\mathcal{G}$ and the results of section 2.4, one
can claim that each data admits a non positive representative. And
then it suffices to show that
$\mathcal{U}=\Theta(\mathcal{F},\mathcal{G})=cl(\theta_{\veps}(f_{\veps},g_{\veps}))_{\veps}$
admits a non positive representative, since a non positive
representative of $\mathcal{U}$ is also one for
$\mathcal{I}\left(\mathcal{U}\right)$. Let  $f_\veps$ and
$g_\veps$ be the non positive representatives of $\mathcal{F}$ and
$\mathcal{G}$, and $u_\veps=\theta_{\veps}(f_{\veps},g_{\veps})$.
Using the maximum's principle as in the proof of proposition 8
with
$$M_\veps=\max\{\sup_{\Omega}f_\veps~,~\sup_{\partial\Omega}g_\veps\}=0,$$
we obtain that $u_\veps \le 0\;a.e.\;\Omega$, and for all $\veps$.
\eproof

\begin{remark}
In fact, we solved the following obstacle problem:\\

For
$\mathcal{F}$$=cl(f_{\veps})_{\veps}\in\mathcal{A}(\mathcal{C},L^{\infty}(\Omega))$
and
$\mathcal{G}$$=cl(g_{\veps})_{\veps}\in\mathcal{A}(\mathcal{C},L^{\infty}(\partial\Omega))$
non positive, find
$\mathcal{U}\in\mathcal{A}(\mathcal{C},H^1(\Omega)\cap
L^{\infty}(\Omega))$ so that
\begin{equation}\label{pbgenobst}
\left\{
\begin{array}{rcl}
-\Delta \Phi(\mathcal{U})+\mathcal{U}&\overset{2}{\simeq}&\mathcal{F} \mbox{ in } \mathcal{A}(\mathcal{C},L^{\infty}(\Om)),\\
\Gamma(\mathcal{U})&=&\mathcal{G} \mbox{ in }
\mathcal{A}(\mathcal{C},L^{\infty}(\partial\Om)),\\
\mathcal{U}&\le&0 \mbox{ in }
\mathcal{A}(\mathcal{C},L^{\infty}(\Om)),
\end{array}
\right.
\end{equation}
which is a generalized version of this one:\\

Find $u:\Omega\mapsto\R$ so that
$$
\left\{
\begin{array}{rcl}
-\Delta \Phi\left(u\right)+u=f & in &\Omega,\\
u=g &  on&\partial\Omega,\\
u\le 0 & on & \Omega,
\end{array}
\right.
$$
where $f$ and $g$ are non positive given functions.
\end{remark}

\end{document}